%% file: StrongGamma_v2.tex
\documentclass[12pt]{amsart}
\usepackage{amsmath, amssymb, latexsym, amsthm, mathrsfs, color, mathtools, tikz,pgfplots,tikz-cd, graphicx, stackengine,url}
\usepackage[top=2.5cm, bottom=2.5cm, left=2.5cm, right=2.5cm]{geometry}

\input{shortcuts}
\usepackage[all]{xy}

\title{}

\author[P.~Szewczak]{Piotr Szewczak}
\address{Piotr Szewczak, Institute of Mathematics, Faculty of Mathematics and Natural Science College of Sciences, Cardinal Stefan Wyszy\'nski University in Warsaw, W\'oycickiego 1$\slash$3, 01--938 Warsaw, Poland
}
\email{p.szewczak@wp.pl}
\urladdr{http://piotrszewczak.pl}

\author[T.~Weiss]{Tomasz Weiss}
\address{Tomasz Weiss, Institute of Mathematics, Faculty of Mathematics and Natural Science College of Sciences, Cardinal Stefan Wyszy\'nski University in Warsaw, W\'oycickiego 1$\slash$3, 01--938 Warsaw, Poland
}
\email{tomaszweiss@o2.pl}

		

\title{Null sets and combinatorial covering properties}

\begin{document}

\maketitle

\begin{abstract}
A subset of the Cantor cube is null-additive if its algebraic sum with any null set is null.
We construct a set of cardinality continuum such that: all continuous images of the set into the Cantor cube are null-additive, it contains a homeomorphic copy of a set that is not null-additive, and it has the property $\gamma$, a strong combinatorial covering property. 
We also construct a nontrivial subset of the Cantor cube with the property $\gamma$ that is not null additive.
Set-theoretic assumptions used in our constructions are far milder than used earlier by Galvin--Miller and Bartoszy\'nski--Rec{\l}aw, to obtain sets with analogous properties.
We also consider products of Sierpi\'nski sets in the context of combinatorial covering properties.
\end{abstract}

\section{Introduction}

Let $\bbN$ be the set of natural numbers and $\PN$ be the power set of $\bbN$.
We identify each set in $\PN$ with its characteristic function, an element of the Cantor cube $\Cantor$; in that way we introduce topology in $\PN$.
The Cantor space $\PN$ with the symmetric difference operation $\oplus$ is a topological group; this operation coincides with the addition modulo $2$ in $\Cantor$.
A set $X\sub\PN$ is \emph{null-additive} in $\PN$ if for any null set $Y\sub \PN$ the set $X\oplus Y:=\sset{x\oplus y}{x\in X, y\in Y}$ is null.
In an analogous way, define null-additive subsets of the real line with the addition $+$ as a group operation.
As we see in the forthcoming Theorem~\ref{thm:nulladd}, it is relatively consistent with ZFC that null-additive subsets of $\PN$ are not preserved by homeomorphisms into $\PN$.
Subsets of the real line whose all continuous images into the real line are null-additive were considered by Galvin and Miller~\cite{gami}; to this end they used combinatorial covering properties.

By \emph{space} we mean a Tychonoff topological space.
A \emph{cover} of a space is a family of proper subsets of the space whose union is the entire space.
An \emph{open} cover of a space is a cover whose members are open subsets of the space.
A cover of a space is an \emph{$\w$-cover} if each finite subset of the space is contained in a set from the cover and it is a \emph{$\gamma$-cover} if it is infinite and each point of the space belongs to all but finitely many sets from the cover. 
A space has the property $\gamma$ if every open $\w$-cover of the space contains a $\gamma$-cover.
This property was introduced by Gerlits and Nagy in the context of local properties of functions spaces~\cite{gn}.
They proved that a space $X$ has the property $\gamma$ if and only if the space $\Cp(X)$ of all continuous real-valued functions defined on $X$ with the pointwise convergence topology is \emph{Fr\'echet--Urysohn}, i.e., each point in the closure of a subset of $\Cp(X)$ is a limit of a sequence from the set~\cite[Theorem~2]{gn}.
Galvin and Miller observed that for a subset of the real line $X$ with the property $\gamma$ and a meager subset of the real line $Y$, the set $X + Y:=\sset{x+y}{x\in X, y\in Y}$ is meager.
They pointed out that they were unable to prove an analogous statement for null sets, and thus they introduced a formally stronger property than $\gamma$~\cite[p. 152]{gami}.
For a natural number $n$, an open cover of a space is an \emph{$n$-cover} if each $n$-elements subset of the space is contained in a member of the cover.
A space $X$ has the property \emph{strongly $\gamma$} if there is an increasing sequence $f\in\NN$ such that for each sequence $\eseq{\cU}$ where $\cU_n$ is an $f(n)$-cover of $X$, there are sets $\seleseq{U}$ such that the family $\sset{U_n}{n\in\bbN}$ is a $\gamma$-cover of $X$.

Galvin and Miller proved that any subset of the real line with the property strongly $\gamma$ is null-additive.
The property strongly $\gamma$ is preserved by continuous mappings, and thus any continuous image of a set with the property strongly $\gamma$, into the real line, is null-additive~\cite[Theorem~7]{gami}.
Under Martin Axiom, Galvin and Miller constructed a  subset of $\PN$ of cardinality continuum with the property strong-$\gamma$~\cite[Theorem~8]{gami}.
In 1996, under some set-theoretic assumption, Bartoszy\'nski and Rec{\l}aw~\cite{br} constructed a subset of $\PN$ of cardinality continuum with the property $\gamma$  that is not null-additive (in particular, they separated the properties $\gamma$ and strongly $\gamma$).

Assuming only an equality between some cardinal characteristics of the continuum, we construct a subset of $\PN$ of cardinality continuum with the property $\gamma$ whose all continuous images into $\PN$ are null-additive and it contains a homeomorphic copy of a set that is not null-additive.  
We also weaken a set-theoretic assumption in the result of Bartoszy\'nski and Rec{\l}aw.
In the both cases we use combinatorial methods of construction of subsets of $\PN$ with the property $\gamma$ invented by Tsaban (\cite[Theorem~3.6]{ot},~\cite[Theorem~6]{sss}) and developed by W{\l}udecka and the first named author~\cite{szw}.
We also use a combinatorial covering characterization of null-additive subsets of $\PN$ given by Zindulka~\cite{zind}.

\section{Null-additive sets with the property $\gamma$}
\label{sec:null-add}

Let $\roth$ be the family of infinite subsets of $\bbN$.
Each set in $\roth$ we identify with an increasing function from $\NN$.
Depending on the interpretation, points of $\roth$ are referred to as sets
or functions.
For natural numbers $n,m$ with $n<m$, define $[n,m):=\sset{i\in\bbN}{n\leq i<m}$.
For sets $a$ and $b$, we write $a\as b$, if the set $a\sm b$ is finite.
A \emph{pseudointersection} of a family of infinite sets is an infinite set $a$ with $a\as b$ for all sets $b$ in the family. A family of infinite sets is \emph{centered} if the finite intersections of its elements, are infinite. 
Let $\fp$ be the minimal cardinality of a subfamily of $\roth$ that is centered and has no pseudointersection.
Let $\Fin$ be the family of finite subsets of $\bbN$.
The following notion plays a crucial role in our constructions.

\bdfn[{Szewczak, W{\l}udecka~\cite{szw}}]
A set $X\sub \roth$ with $\card{X}\geq \fp$ is a \emph{$\fp$-generalized tower} if for each function $a\in\roth$, there are sets $b\in\roth$ and $S\sub X$ with $\card{S}<\fp$ such that 
\[
x\cap \Un_{n\in b}[a(n),a(n+1))\in\Fin
\]
for all sets $x\in X\sm S$.
\edfn

For functions $f,g\in\roth$, we write $f\les g$, if the set $\sset{n}{f(n)>g(n)}$ is finite.
A subset of $\roth$ is \emph{unbounded} if for any function $g\in\roth$, there is a function $f$ in the set with $f\not\les g$.
Let $\fb$ be the minimal cardinality of an unbounded subset of $\roth$.
The existence of a $\fp$-generalized tower in $\roth$ is independent of ZFC, i.e., it is equivalent to the equality $\fp=\fb$~\cite[Lemma~3.3]{ot}.
Let $\nonNadd$ be the minimal cardinality of a subset of $\PN$ that is not null-additive.

\bthm\label{thm:nulladd}
Assume that $\fp=\non(\nadd)=\fc$.
There is a set $X\sub\roth$ such that
\be
\item The set $X$ is a $\fp$-generalized tower.
\item The set $X\cup\Fin$ has the property $\gamma$.
\item All continuous images of the set $X\cup\Fin$ into $\PN$ are null-additive.
\item The set $X$ is homeomorphic to a subset of $\PN$ that is not null-additive.
\ee
\ethm

We need the following notions and auxiliary results.

For a class $\cA$ of covers of spaces and a space $X$, let $\cA(X)$ be the family of all covers of $X$ from the class $\cA$.
Let $\eseq{\cA}$ and $\cB$ be classes of covers of spaces.
A space $X$ satisfies
\[
\sone(\{\cA_n\}_{n\in\bbN}, \cB)
\]
if for each sequence $\eseq{\cU}$ with $\cU_n\in\cA_n(X)$ for each $n$, there are sets $\seleseq{U}$ such that $\sset{U_n}{n\in\bbN}\in\cB(X)$.
Let $\Op_n$ be the class of all open $n$-covers of spaces for each $n$ and $\Ga$ be the class of all open $\gamma$-covers of spaces.
Using the above notion, the strongly $\gamma$ property is the property $\strongga$.
Let $d$ be a metric in $\PN$ that coincides with the standard metric in $\Cantor$, that is, for points $x,y\in \PN$, let
\[
d(x,y):=
\begin{cases}
2^{-\min((x\cup y)\sm (x\cap y))}, &\text{ if }x\neq y,\\
0, &\text{ if }x=y.
\end{cases}
\]
Let $\epsilon$ be a positive real number.
For a point $x\in \PN$, let $B(x,\epsilon):=\sset{y\in \PN }{d(x,y)<\epsilon}$ and for a set $F\sub \PN$, let $B(F,\epsilon):=\Un\sset{B(x,\epsilon)}{x\in F}$.
Fix a natural number $n$.
An open cover of a subspace of $\PN$ is a \emph{uniform $n$-cover} if there is a positive real number $\epsilon$ such that for each $n$-elements subset $F$ of the space, the set $B(F,\epsilon)$ is contained in a member of the cover.
Let $\Op_n^{\mathsf{unif}}$ be the family of all uniform $n$-covers of subspaces of $\PN$.
In 1995, Shelah characterized null-additive subsets of $\PN$~\cite{shelah}.
Using the Shelah result, Zindulka characterized null-additive subsets of $\PN$, in terms of  combinatorial covering properties.

\bthm[Zindulka~\cite{zind}]\label{thm:zind}
A subset of $\PN$ is null-additive if and only if it satisfies $\sone(\{\Op_n^{\mathsf{unif}}\}_{n\in\bbN},\Ga)$.
\ethm

For a set $X$ and a natural number $n$, let $\Pof_n(X)$ be the family of all $n$-elements subsets of $X$.
By a similar argument as in the proof of a result of Tsaban~\cite[Theorem~2.1]{Tsstrongga}, we have the following result.

\blem\label{lem:zind}
For a set $X\sub\PN$ the following assertions are equivalent:

\be
\item The set $X$ satisfies $\nullzind$.
\item There is a function $f\in\roth$ such that the set $X$ satisfies $\nullzindf{f(n)}$.
\item For each function $f\in\roth$, the set $X$ satisfies $\nullzindf{f(n)}$.\qed
\ee
\elem

\bprp\label{prp:strongNadd}
Let $X\sub\PN$ be a set and assume that there is a function $f\in\roth$ such that for each sequence $\eseq{\cU}$ with $\cU_n\in\Op_{f(n)}(X)$ for each $n$, there are sets $\seleseq{U}$ and a set $X'\sub X$ with $\card{X'}<\nonNadd$ such that $\sset{U_n}{n\in\bbN}\in\Ga(X\sm X')$.
Then all continuous images of the set $X$ into $\PN$ are null-additive.
\eprp

\bpf
Let $\eseq{\cU}$ be a sequence such that $\cU_n\in \Op^{\mathsf{unif}}_{2f(n)}(X)$ for each $n$.
Fix a natural number $n$.
We may assume that there is a positive real number $\epsilon_n$ such that
\[
\cU_n=\sset{B(A,\epsilon_n)}{A\in\Pof_{2f(n)}(X)}.
\]
We have $\cU'_n:=\sset{B(F,\epsilon_n)}{F\in\Pof_{f(n)}(X)}\in \Op^{\mathsf{unif}}_{f(n)}(X)$.
By the assumption, there are sets $F_1\in\Pof_{f(1)}(X), 
F_2\in\Pof_{f(2)}(X), \dotsc$ and a set $X'\sub X$ with $\card{X'}<\nonNadd$ such that 
\[
\sset{B(F_n,\epsilon_n)}{n\in\bbN}\in\Ga(X\sm X').
\]
Since $\card{X'}<\nonNadd$, there are sets $F_1'\in\Pof_{f(1)}(X'), 
F_2'\in\Pof_{f(2)}(X'), \dotsc$ such that 
\[
\sset{B(F_n',\epsilon_n)}{n\in\bbN}\in\Ga(X').
\]
For each natural number $n$, there is a set $A_n\in\Pof_{2f(n)}(X)$ such that $F_n\cup F_n'\sub A_n$ and
\[
B(F_n,\epsilon_n)\cup B(F_n',\epsilon_n)\sub B(A_n,\epsilon_n)\in\cU_n.
\]
We have
\[
\sset{B(A_n,\epsilon_n)}{n\in\bbN}\in \Ga(X)
\]
By Lemma~\ref{lem:zind} and Theorem~\ref{thm:zind}, the set $X$ is null-additive.
\epf

By Theorem~\ref{thm:zind}, we have the following result.

\bprp\label{prp:unif}
Null-additivity in $\PN$ is preserved by uniformly continuous functions. \qed
\eprp

\blem\label{lem:n+k}
Let $X$ be a space and $n,k$ be natural numbers.
If $\cU\in\Op_{n+k}(X)$, then each $n$-elements subset of $X$ is contained in at least $k$ pairwise different sets from $\cU$.\qed
\elem

For the remaining part of this section, let $f\in\roth$ be a sequence such that $f(n+1)>2^{f(n)+2n}+n$ for all natural numbers $n$. For a set $d\in\roth$, let $\rothx{d}$ be the family of all infinite subsets of $d$.

\blem\label{lem:GMstrong}
Let $\eseq{\cU}$ be a sequence of families of open sets in $\PN$ such that  $\cU_n\in\Op_{f(n)}(\Fin)$ for all natural numbers $n$ and $d\in\roth$.
There are an element $x\in\rothx{d}$ and pairwise different sets $\seleseq{U}$ such that for each natural number $n$ and each set $y\in\PN$:
\[
\text{If }y\cap [f(n),x(n+1))\sub \{x(1),x(1)+1,\dotsc,x(n), x(n)+1\},\text{ then }y\in U_{n+1}.
\]
In particular, 
\[
\sset{U_n}{n\in\bbN}\in\Ga(\sset{y\in\PN}{y\as x\cup(x+1)}).
\]
\elem

\bpf
Let $x(1):=a(1)$ and $U_1\in\cU_1$.
Fix a natural number $n$ and assume that natural numbers $x(1),\dotsc,x(n)$ and sets $U_1\in\cU_1,\dotsc,U_n\in\cU_n$ have already been defined.
We have
\[
\card{\Pof([1,f(n))\cup \{x(1),x(1)+1,\dotsc,x(n), x(n)+1\})}\leq 2^{f(n)+2n}
\]
and $\cU_{n+1}\in \Op_{f(n+1)}(\Fin)$.
By Lemma~\ref{lem:n+k} there is a set $U_{n+1}\in\cU_{n+1}\sm\{U_1,\dotsc, U_n\}$ such that 
\[
\Pof([1,f(n))\cup \{x(1),x(1)+1,\dotsc,x(n), x(n)+1\})\sub U_{n+1}.
\]
For each set 
\[
s\in \Pof([1,f(n))\cup \{x(1),x(1)+1,\dotsc,x(n), x(n)+1\}),
\]
there is a natural number $m_s\in a$ with $m_s>\max\{f(n),x(n)+1\}$ such that for each set $y\in\PN$:
\[
\text{If }y\cap [1,m_s)=s,\text{ then }y\in U_{n+1}.
\]
Define
\[
x(n+1):=\max\sset{m_s}{s\in  \Pof([1,f(n))\cup \{x(1),x(1)+1,\dotsc,x(n), x(n)+1\})}.
\]

Fix a natural number $n$.
Let $y\in\PN$ be a set such that
\[
y\cap [f(n), x(n+1))\sub \{x(1),x(1)+1,\dotsc,x(n), x(n)+1\}.
\]
Then the set $s:=y\cap [1, x(n+1))$ belongs to $\Pof([1,f(n))\cup \{x(1),x(1)+1,\dotsc,x(n), x(n)+1\})$, and
\[
y\cap [1, m_s)=y\cap[1,x(n+1))=s.
\]
Thus, $y\in U_{n+1}$.
\epf

\blem[{Folklore~\cite[Lemma~2.13]{MHP}}]\label{lem:unbdd}
A set $X\sub\roth$ is unbounded if and only if for each function $a\in\roth$, there are a set $b\in\roth$ and an element $x\in X$ such that 
\[
x\cap \Un_{n\in b}[a(n),a(n+1))=\emptyset
\]
for all $x\in X$.
\elem

For elements $b\in\roth$ and $c\in\Cantor$, let $2b$, $b+1$ and $b+c$ be elements in $\roth$ such that $(2b)(n):=2b(n)$, $(b+1)(n):=b(n)+1$ and $(b+c)(n):=b(n)+c(n)$ for all natural numbers $n$.

\begin{proof}[{Proof of Theorem~\ref{thm:nulladd}}]
Construction of a set $X$:
Let $\sset{(\eseq{\cU^{(\alpha)}})}{\alpha<\fc}$ be an enumeration of all sequences of families of open sets in $\PN$ such that $\cU^{(\alpha)}_n\in\Op_{f(n)}(\Fin)$ for all ordinal numbers $\alpha<\fc$ and natural numbers $n$.
Let $\roth=\sset{d_\alpha}{\alpha<\fc}$.

Apply Lemma~\ref{lem:GMstrong} to the sequence $\eseq{\cU^{(0)}}$ and to a set $d\in\rothx{2\bbN}$ with $d_0\les d$.
Then there are sets $U^{(0)}_1\in\cU^{(0)}_1,U^{(0)}_2\in\cU^{(0)}_2,\dotsc$ and an element $\pick_0\in\rothx{2\bbN}$ such that 
\[
\sset{U_n^{(0)}}{n\in\bbN}\in\Ga(\sset{\gen\in\PN}{\gen\as \pick_0\cup(\pick_0+1)}).
\]
Fix an ordinal number $\alpha<\fc$.
Assume that elements $\pick_\beta\in\rothx{2\bbN}$ with $d_\beta\les \pick_\beta$ and sets $U^{(\beta)}_1\in\cU^{(\beta)}_1, U^{(\beta)}_2\in\cU^{(\beta)}_2,\dotsc$ with
\[
\sset{U^{(\beta)}_n}{n\in\bbN}\in\Ga(\sset{\gen\in\PN}{\gen\as\pick_\beta\cup(\pick_\beta+1)}
\]
have already been defined for all ordinal numbers $\beta<\alpha$ such that for $\beta,\beta'<\alpha$ if $\beta<\beta'$, then $\pick_\beta\aspst \pick_{\beta'}$.
Let $d\in\rothx{2\bbN}$ be a pseudointersection of the family $\sset{\pick_\beta}{\beta<\alpha}$ with $d_\alpha\les d$.
Apply Lemma~\ref{lem:GMstrong} to the sequence $\eseq{\cU^{(\alpha)}}$ and to the element $d$.
Then there are an element $\pick_\alpha\in\rothx{d}$ and sets $U^{(\alpha)}_1\in\cU^{(\alpha)}_1, U^{(\alpha)}_2\in\cU^{(\alpha)}_2,\dotsc$  such that
\begin{equation}
\label{eq:ind}
\sset{U_n^\alpha}{n\in\bbN}\in\Ga(\sset{y\in\PN}{y\as x_\alpha\cup(x_\alpha+1)}).
\end{equation}

Let $\Cantor=\sset{c_\alpha}{\alpha<\fc}$.
Define
\[
X:=\sset{\pick_\alpha+c_\alpha}{\alpha<\fc}
\] 
By the construction, for all ordinal numbers $\alpha,\beta$ with $\alpha\leq\beta<\fc$, we have
\begin{equation}\label{eq:as}
\pick_\beta+c_\beta\sub \pick_\beta\cup (\pick_\beta+1)\as \pick_\alpha\cup(\pick_\alpha+1).
\end{equation}

(1)
Let $a\in \roth$.
Since the set $\sset{\pick_\alpha}{\alpha<\fc}$ is unbounded, the set $\sset{\pick_\alpha\cup(\pick_\alpha+1)}{\alpha<\fc}$ is unbounded, too.
By Lemma~\ref{lem:unbdd}, there are an ordinal number $\alpha<\fc$ and a set $b\in\roth$ such that 
\[
(\pick_\alpha\cup (\pick_\alpha+1))\cap\Un_{n\in b}[a(n),a(n+1))=\emptyset.
\]
By~\eqref{eq:as}, we have
\[
(\pick_\beta+c_\beta) \cap\Un_{n\in b}[a(n),a(n+1))\as (\pick_\alpha\cup (\pick_\alpha+1))\cap\Un_{n\in b}[a(n),a(n+1))=\emptyset
\]
for all ordinal numbers $\beta$ with $\alpha\leq\beta<\fc$.
Thus, the set $X$ is a $\fp$-generalized tower.

(2) By~(1), the set $X$ is a  $\fp$-generalized tower, and thus the set $X\cup\Fin$ has the property $\gamma$~\cite[Theorem~4.1(1)]{szw}.

(3)
The set $X\cup\Fin$ satisfies the property from Proposition~\ref{prp:strongNadd}:
Let $\eseq{\cU}$ be a sequence of open families in $\PN$ such that $\cU_n\in\Op_{f(n)}(X\cup\Fin)$ for all natural numbers $n$.
There is an ordinal number $\alpha<\fc$ such that this sequence is equal to the sequence $\eseq{\cU^{(\alpha)}}$.
By~\eqref{eq:ind} and~\eqref{eq:as}, we have
\[
\sset{U^{(\alpha)}_n}{n\in\bbN}\in\Ga(\sset{\pick_\beta+c_\beta}{\alpha\leq\beta<\fc}).
\]

(4)
In our proof we use Rothberger's trick~\cite[Theorem~9.4]{mill}.
The map $\Phi\colon \PN\to\Cantor$ unifying each subset of $\bbN$ with its characteristic function is a homeomorphism and an isometry with respect to the standard metric in $\Cantor$ and the metric $d$ on $\PN$ defined in Section~\ref{sec:null-add}.
The set 
\[
G:=\smallmedset{((\Phi(x))(1), x(1)_{\bmod 2},(\Phi(x))(2), x(2)_{\bmod 2},\dotsc)}{x\in X},
\]
is a homeomorphic copy of the set $X$.
Since the projection $\pi\colon G\to \{0,1\}^{2\bbN}$ is onto and uniformly continuous, the map
\[
\Phi\inv\circ \pi \circ \Phi\colon \Phi\inv[G]\to \PN
\]
is onto and uniformly continuous, too.
By Proposition~\ref{prp:unif}, the set $\Phi\inv[G]$ is not null additive.
The set $\Phi\inv[G]$ is homeomorphic to $G$.
\epf

\section{A $\gamma$-set that is not null-additive}

Let $\non(\cN)$ be the minimal cardinality of a subset of $\PN$ that is not null.
A set $\sset{x_\alpha}{\alpha<\fp}\sub\roth$ is a $\fp$-unbounded tower if it is unbounded and for all ordinal numbers $\alpha,\beta<\fp$ with $\alpha<\beta$, we have $x_\alpha\aspst x_\beta$.
A $\fp$-unbounded tower exists if and only if $\fp=\fb$~\cite[Lemma~3.3]{ot}.

\setcounter{equation}{0}

\bthm\label{thm:tower_non_Nadd}
Assume that $\fp=\fb=\non(\cN)$.
There is a $\fp$-generalized tower $X\sub\roth$ that is not null-additive.
In particular, the set $X\cup\Fin$ is a nontrivial set with the property $\gamma$ that is not null-additive.
\ethm

\bpf
Let $f\in\roth$ be a function such that $f(n):=\sum_{i=0}^{n}2^i$ for all natural numbers $n$. 
Define 
\[
Y_n:=\sset{x\in\PN}{x\cap [f(n),f(n+1))=\emptyset}
\]
for all natural numbers $n$.
Each set $Y_n$ is clopen and has measure less or equal than $2^{-(n+1)}$, and thus the set $Y:=\bigcap_m\bigcup_{n>m}Y_n$ is null.
Let $T=\sset{t_\alpha}{\alpha<\fp}$ be a $\fp$-unbounded tower in $\roth$ and $Z=\sset{z_\alpha}{\alpha<\fp}$ be a nonnull set in $\PN$.

Define
\[
x_\alpha:=z_\alpha\cap \Un_{n\in\bbN}[f(t_\alpha(2n)),f(t_\alpha(2n)+1)
\cup\sset{f(t_\alpha(2n+1))}{n\in\bbN}
\]
for all ordinal numbers $\alpha<\fp$ and $X:=\sset{x_\alpha}{\alpha<\fp}$.

The set $X$ is a $\fp$-generalized tower:
For each ordinal number $\alpha$ with $\alpha<\fp$, we have
\[
\sset{f(t_\alpha(2n+1))}{n\in\bbN}\sub x_\alpha,
\]
and thus $X\sub\roth$.
Let $a\in\roth$.
There is a set $c\in\roth$ such that
\[
\card{[f(c(k)),f(c(k+1)))\cap a}\geq 2
\]
for all natural numbers $k$.
By Lemma~\ref{lem:unbdd}, there are a set $d\in\roth$ and an ordinal number $\alpha<\fp$ such that 
\begin{equation}\label{eq:tower}
t_\alpha\cap\Un_{k\in d}[c(k),c(k+1))=\emptyset.
\end{equation}
Then there is a set $b\in\roth$ such that 
\[
\Un_{n\in b}[a(n),a(n+1))\sub \Un_{k\in d}[f(c(k)),f(c(k+1))).
\]
Fix an ordinal number $\beta$ with $\alpha\leq\beta<\fp$.
We have
\begin{gather*}
x_\beta\cap \Un_{n\in b}[a(n),a(n+1))\sub\\
\Un_{n\in t_\beta}[f(n),f(n+1))\cap \Un_{n\in b}[a(n),a(n+1))\as
\Un_{n\in t_\alpha}[f(n),f(n+1))\cap \Un_{k\in d}[f(c(k)),f(c(k+1))).
\end{gather*}
By~\eqref{eq:tower}, the latter intersection is empty.

We have $Z\sub X\oplus Y$:
Fix an ordinal number $\alpha<\fp$.
Since 
\[
x_\alpha\cap \Un_{n\in\bbN}[f(t_\alpha(2n)),f(t_\alpha(2n)+1)= z_\alpha\cap \Un_{n\in\bbN}[f(t_\alpha(2n)),f(t_\alpha(2n)+1),
\]
we have
\[
(x_\alpha\oplus z_\alpha)\cap \Un_{n\in\bbN}[f(t_\alpha(2n)),f(t_\alpha(2n)+1)=\emptyset ,
\]
and thus $x_\alpha\oplus z_\alpha\in Y$.
Then there is an element $y_\alpha\in Y$ such that $x_\alpha\oplus z_\alpha=y_\alpha$.
Thus, $z_\alpha=x_\alpha\oplus y_\alpha$.

Since $X$ is a $\fp$-generalized tower, the set $X\cup\Fin$ has the property $\gamma$~\cite[Theorem~4.1(1)]{szw}.
Since $Z\sub (X\cup\Fin)\oplus Y$, the set $X\cup\Fin$ is not null-additive.
\epf

\brem
The assumption of Theorem~\ref{thm:tower_non_Nadd} is valid assuming \CH{} or Martin Axiom and in the following model.
Let $\bbP_{\aleph_1}$ be an $\aleph_1$-iteration with finite support of a measure algebra and $G$ be a generic filter in $\bbP_{\aleph_1}$.
Let $M$ be a model of ZFC and $\fc=\aleph_2$.
In the model $M[G]$, the assumption $\fp=\fb=\non(\cN)$ from Theorem~\ref{thm:tower_non_Nadd} is true:
Let $\non(\cM)$ be the minimal cardinality of a subset of $\PN$ that is not meager.
Adding $\aleph_1$ Cohen reals, we have $\non(\cM)=\aleph_1$ and adding $\aleph_1$ Solovay reals, we have $\non(\cN)=\aleph_1$.
Thus, $\fp=\aleph_1$.
Since $\non(\cM)=\aleph_1$, we have $\fb=\aleph_1$.
\erem

\bthm
Assume that $\fp=\fb=\non(\cN)$.
There is a nontrivial subset of the real line with the property $\gamma$ that is not null-additive.
\ethm

\bpf
We use notions from the paper of the second named author~\cite[Corollary~11]{addWeiss}.
Let $X\sub \roth$ be a set from Theorem~\ref{thm:tower_non_Nadd}.
Let $p\colon \PN\to A$ be a homeomorphism. 
Since the map $f\circ g\inv\circ p$ is continuous, the image $Y$ of the set $X\cup\Fin$ under the map $f\circ g\inv\circ p$ is a subset of the real line with the property $\gamma$.
Suppose that the set $Y$ is null-additive.
Then the set $g\circ f\inv[Y]$ is null-additive~\cite[Theorem~12]{addWeiss}.
Since the map $p\inv$ is unifromly continuous and the set $X\cup\Fin$ is not null-additive, the set $p[X]=(g\circ f\inv[Y])$ is not null-additive, too, a contradiction.
\epf

\bthm
It is consistent with ZFC that there is a $\fp$-unbounded tower in $\roth$ and $\nonNadd<\fp$.
\ethm

\bpf
The dual Borel conjecture is the statement that for any set $X\sub\PN$ with cardinality $\aleph_1$, there is a null set $N\sub \PN$ such that $X\oplus N=\PN$.
By the result of Judah and Shelah~\cite[Theorem~8.5.23]{BJ} there is a model for ZFC satisfying Martin Axiom for $\sigma$-centered sets and the dual Borel conjecture.
In that model, we have $\fp=\aleph_2=\fb$, an thus a $\fp$-unbounded tower exists. On the other hand, since the dual Borel conjecture holds, we have $\nonNadd=\aleph_1$.

\epf

\section{Nonproductivity of Sierpi\'nski-type sets}

Let $\Op$ be the class of open covers of spaces.
A space $X$ satisfies Menger's property $\men$ if for each sequence $\eseq{\cU}\in\Op(X)$, there are finite sets $\cF_1\sub\cU_1,\cF_2\sub \cU_2,\dotsc$ such that $\Un_n\cF_n\in\Op(X)$.
A space $X$ satisfies Hurewicz's property $\hur$ if for each sequence $\eseq{\cU}\in\Op(X)$, there are finite sets $\cF_1\sub\cU_1, \cF_2\sub\cU_2,\dotsc$ such that $\sset{\Un\cF_n}{n\in\bbN}\in \Ga(X)$.
The property $\hur$ implies $\men$ and it generalizes $\sigma$-compactness.
An uncountable subset of $\PN$ is a \emph{Sierpi\'nski} set if its intersection with any null set is at most countable; its existence is independent of ZFC.
Any Sierpi\'nski set satisfies $\hur$ but it is not $\sigma$-compact.
Assuming that \CH{} holds, there are two Sierpi\'nski sets whose product space does not even satisfy $\men$~\cite[p. 250]{coc2}.

A category theoretic counterpart to a Sierpi\'nski set is a \emph{Luzin set}, i.e.,  an uncountable subset of $\PN$ whose intersection with any meager set is at most countable. Each Luzin set satisfies $\men$ but no $\hur$.
A \emph{set of reals} is a space homeomorphic to a subset of the real line.
Let $\bfP$ be a property of spaces.
A space is \emph{productively $\bfP$} if its product space with any space satisfying $\bfP$, satisfies $\bfP$.
Assuming $\aleph_1=\cf(\fd)$, in the class of sets of reals, no Luzin set is productively $\men$~\cite[Corollary~2.11.]{ST}.
There are open problems~(\cite[Problem~7.5]{ST},~\cite[Problem~5.5]{pMGen}), whether in the class of sets of reals (or in the class of general topological spaces) for any Sierpi\'nski set $X$, there is a space $Y$ satisfying $\hur$ such that the product space $X\x Y$ does not satisfy $\hur$, and what if we assume \CH{}?
We consider an analogous problem with respect to combinatorial covering properties, stronger than $\hur$.

Let $\cA$, $\cB$ be classes of covers of spaces.
A space $X$ satisfies $\sone(\cA,\cB)$ if for each sequence $\eseq{\cU}\in\cA(X)$, there are sets $\seleseq{U}$ such that $\sset{U_n}{n\in\bbN}\in\cB(X)$.
Let $\Om$ be the class of all open $\w$-covers of spaces.
A \emph{Borel} cover of a space is a cover whose members are Borel subsets of the space.
Let $\Gab$ and $\Omb$ be classes of all  countable Borel $\gamma$-covers and countable Borel $\w$-covers of spaces, respectively. 
The property $\gamma$, considered in the previous sections, is one of the classic properties in the \emph{selection principles} theory, it is equivalent to the property $\swg$~\cite[Theorem~2]{gn}.
A set of reals is \emph{totally imperfect} if it does not contain an uncountable compact set and it is \emph{perfectly meager} if its intersection with a perfect set is meager in this perfect set.
The following diagram presents relations between considered properties~\cite{coc1},~\cite{coc2}.

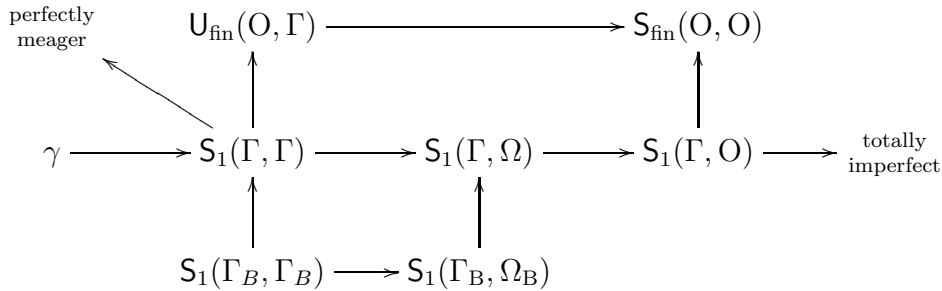
\[
\xymatrix{
	\genfrac{}{}{0pt}{}{\text{perfectly}}{\text{meager}}&\hur \ar[rr]& &\men&\\
	\gamma \ar[r]&\sgg\ar[u]\ar[r]\ar[lu] &\sgw\ar[r] &\sgo \ar[u]\ar[r]& 	\genfrac{}{}{0pt}{}{\text{totally}}{\text{imperfect}}\\
	& \sbgg\ar[u]\ar[r] &\sbgw\ar[u]&&
}
\]

Let $\cov(\cN)$ be the minimal cardinality of a family of null subsets of $\PN$ whose union is $\PN$ and $\cof(\cN)$ be the minimal cardinality of a family of null subsets of $\PN$ such that any null subset of $\PN$ is contained in a set from the family.
For an uncountable ordinal number $\kappa$, a set $X\sub\PN$ is a \emph{$\kappa$-Sierpi\'nski set} if $\card{X}\geq\kappa$ and for any null set $Y\sub\PN$, we have $\card{X\cap Y}<\kappa$.
Every $\fb$-Sierpi\'nski set satisfies $\sbgg$~(\cite[Theorem~2.9]{coc2},~\cite[Theorem~2.4]{MHP}).

\bthm\label{thm:sierp}\mbox{}
\be
\item 
Assume that $\cov(\cN)=\cof(\cN)=\fb$.
For every $\fb$-Sierpi\'nski set $S$, there is a $\fb$-Sierpi\'nski set $S'$ such that the product space $S\x S'$ does not satisfy $\sgg$.
\item Assume that $\cov(\cN)=\cof(\cN)=\fd=\fc$ and the cardinal number $\fc$ is regular.
For every $\fc$-Sierpi\'nski set $S$, there is a $\fc$-Sierpi\'nski set $S'$ such that the product space $S\x S'$ does not satisfy $\sgo$.
\ee
\ethm

In order to prove Theorem~\ref{thm:sierp}, we need the following Lemma.

\blem\label{lem:sierp}
Assume that $\cov(\cN)=\cof(\cN)$ and the cardinal number $\cov(\cN)$ is regular. 
For every $\cov(\cN)$-Sierpi\'nski set $S$ and every set $Y$ of cardinality at most $\cov(\cN)$, there is a $\cov(\cN)$-Sierpi\'nski set $S'$ such that $Y\sub S\oplus S'$.
\elem

\bpf
Let $\sset{N_\alpha}{\alpha<\cov(\cN)}$ be a cofinal family of null sets in $\PN$ and $Y=\sset{y_\alpha}{\alpha<\cov(\cN)}\sub\PN$.
Fix an ordinal number $\alpha<\cov(\cN)$ and assume that elements $x_
\beta\in\PN$ have already been defined for all ordinal numbers $\beta$ with $\beta<\alpha$.
The set $y_\alpha \oplus S$ is a $\cov(\cN)$-Sierpi\'nski.
Since the cardinal number $\cov(\cN)$ is regular, the union $\Un_{\beta<\alpha}N_\beta\cup\sset{x_\beta}{\beta<\alpha}$ cannot cover any $\cov(\cN)$-Sierpi\'nski set.
Then there is an element 
\[
x_\alpha\in (y_\alpha \oplus S)\sm \Un_{\beta<\alpha}N_\beta\cup\sset{x_\beta}{\beta<\alpha}.
\]
By the construction the set $S':=\sset{x_\alpha}{\alpha<\cov(\cN)}$ is a $\cov(\cN)$-Sierpi\'nski set and for each ordinal number $\alpha<\cov(\cN)$, there is an element $s_\alpha\in S$ such that $y_\alpha=x_\alpha\oplus s_\alpha$.
Thus, $Y\sub S\oplus S'$.
\epf

Let $\non(\cM)$ be the minimal cardinality of a nonmeager subset of $\PN$.

\begin{proof}[{Proof of Theorem~\ref{thm:sierp}}]
(1) By the Cicho\'n diagram, we have $\cov(\cN)\leq\non(\cM)\leq\cof(\cN)$, and thus $\non(\cM)=\fb$.
Let $Y\sub\PN$ be a nonmeager set of cardinality $\fb$.
By Lemma~\ref{lem:sierp}, there is a $\fb$-Sierpi\'nski set $S'$ such that $Y\sub S\oplus S'$.
Every set with the property $\sgg$ is meager and the property $\sgg$ is preserved by continuous functions.
Since $S\oplus S'$ is a nonmeager continuous image of the product space $S\x S'$, the product $S\x S'$ does not satisfy $\sgg$.	

(2) 
Every set satisfying $\sgo$ is totally imperfect.
Let $Y=\roth$.
Proceed analogously as in (1).
\epf

\brem
The assumptions of Theorem~\ref{thm:sierp} are valid assuming \CH{} or Martin Axiom.
In a model obtained by $\aleph_2$-iteration of Sacks forcing with countable supports, we have $\cov(\cN)=\non(\cM)=\cof(\cN)=\fb=\aleph_1$.
In a model obtained by $\aleph_2$-iteration of amoeba forcing with finite supports, we have $\cov(\cN)=\cof(\cN)=\fd=\fc$.
\erem

It is not known whether, in the class of totally imperfect sets, the properties $\sgo$ and $\men$ are different. Theorem~\ref{thm:sierp} can be useful to solve this problem.

\brem
Assume that \CH{} holds. If there is a Sierpi\'nski set whose product space with any Sierpi\'nski set satisfies $\men$, then, in the class of totally imperfect sets, the properties $\men$ and $\sgo$ are different.
\erem

For a finite set $F\sub\roth$ and a set $X\sub\roth$ let $\max[F]:=\sset{\max_{f\in F}f(n)}{n\in\bbN}$, an element of $\roth$, and $\maxfin[X]:=\sset{\max[F]}{F\text{ is a finite subset of }X}$.
For elements $f,g\in\roth$ we write $f\leinf g$ if the set $\sset{n}{f(n)\leq g(n)}$ is infinite.

Using similar ideas as in the proof that any Sierpi\'nski set satisfies $\sgg$~\cite[Theorem~2.9, Theorem~2.10]{coc2}, we obtain the following results.

\bprp\label{prp:dsierp}
Every $\fd$-Sierpi\'nski set satisfies $\sbgw$.
\eprp

\bpf
Let $S$ be a $\fd$-Sierpi\'nski set.
There is a positive real number $p$ such that the outer measure of $S$ is equal to $p$.
Let $B$ be a Borel set containing $S$ with $\mu(B)=p$.
Let $\eseq{\cU}\in\borga(S)$ be a sequence of families of Borel subsets of the set $B$.
Let $\cU_n=\sset{U^n_m}{m\in\bbN}$ for all natural numbers $n$, and we may assume that each such a family is increasing (if not consider the family $\sset{\bigcap_{i\geq j}{U^n_i}}{j\in\bbN}$, an increasing countable Borel cover of $S$).
Fix a natural number $k$.
There is a function $f_k\in\NN$ such that $\mu(U^n_{f_k(n)})\geq (1 - \frac{1}{2^{n+k}})p$.
For a set $A_k:=\bigcap_{n}U^n_{f_k(n)}$, we have $\mu(A_k)\geq(1-\frac{1}{2^k})p$.
The set $A=\Un_kA_k$ has measure $p$, and thus $\card{S\sm A}<\fd$.
For each element $x\in S\sm A$, there is a function $f_x\in \NN$ such that $x\in\bigcap_nU^n_{f_x(n)}$.
There is a function $g\in\NN$ such that 
\[
\maxfin[\sset{f_k}{k\in\bbN}\cup\sset{f_x}{x\in S\sm A}]\leinf g.
\]
We have $\sset{U^n_{g(n)}}{n\in\bbN}\in\Om(S)$:
Let $F$ be a finite subset of $S$.
There is a set $a\in\Fin$ such that $F\cap A\sub \Un_{k\in a}A_k$.
Let 
\[
f:=\maxfin[\sset{f_k}{k\in a}\cup\sset{f_x}{x\in S\sm A}].
\]
Then $F\sub U^n_{f(n)}$ for all but finitely many natural numbers $n$.
Since $f\leinf g$, and families $\cU_n$ are increasing, we have $F\sub U^n_{g(n)}$ for infinitely many natural numbers $n$. 
\epf

The following corollary is a straightforward consequence of Theorem~\ref{thm:sierp}.

\bcor\mbox{}
\be
\item 
Assume that $\cov(\cN)=\cof(\cN)=\fb$.
No $\cov(\cN)$-Sierpi\'nski set is productively $\sgg$.
\item Assume that $\cov(\cN)=\cof(\cN)=\fd=\fc$ and the cardinal number $\fc$ is regular.
No $\fc$-Sierpi\'nski set $S$ is productively $\sgw$ or productively $\sgo$.
\ee
\ecor

\section*{Acknowledgments}

We would like to thank the referees for their work on refereeing this paper and their detailed comments and corrections.

\end{document}

%% file: shortcuts.tex
\stackMath
\newcommand\tsup[2][2]{%
 \def\useanchorwidth{T}%
  \ifnum#1>1%
    \stackon[-1pt]{\tsup[\numexpr#1-1\relax]{#2}}{\hspace{1pt}\scriptstyle\sim}%
  \else%
    \stackon[.5pt]{#2}{\hspace{1pt}\scriptstyle\sim}%
  \fi%
}

\newcommand{\nc}{\newcommand}

\newcommand{\swg}{\mathsf{S}_1(\Omega,\Gamma)}
\newcommand{\sgg}{{\mathsf{S}_1(\Ga,\Ga)}}
\newcommand{\strongga}{\sone(\{\Op_n\}_{n\in\bbN}),\Ga)}

\newcommand{\sbgg}{{\mathsf{S}_1(\Ga_B,\Ga_B)}}

\newcommand{\sgo}{{\mathsf{S}_1(\Ga,\Op)}}
\newcommand{\sbgw}{{\mathsf{S}_1(\Ga_{\mathrm{B}},\Om_{\mathrm{B}})}}
\newcommand{\sgw}{{\mathsf{S}_1(\Ga,\Om)}}

\nc{\mc}{\mathcal}
\nc{\thusfar}{\my{--- Edited thus far ---}}
\nc{\lei}{\le^\oo}
\nc{\sqsubs}{\sqsubseteq^*}
\nc{\card}[1]{\left|#1\right|}
\nc{\medcard}[1]{\biggl|\,#1\,\biggr|}
\nc{\smallcard}[1]{|\,#1\,|}
\nc{\bds}{bidirectional $\roth$-scale}
\nc{\bfP}{\mathbf{P}}
\nc{\bfQ}{\mathbf{Q}}
\nc{\bbT}{\mathbb{T}}
\nc{\bbZ}{\mathbb{Z}}
\nc{\bbN}{\mathbb{N}}
\nc{\bbC}{\mathbb{C}}
\nc{\beq}{\begin{equation}}\nc{\eeq}{\end{equation}}
\nc{\mbq}{\mb{?}}
\nc{\mb}[1]{{\mbox{\textbf{#1}}}}
\nc{\nop}{$\times$}
\nc{\fbn}{\!\!\fbox{\!\nop\!}\!\!}
\nc{\yup}{\checkmark}
\nc{\forces}{\Vdash}
\nc{\name}[1]{\dot{#1}}
\nc{\tf}{\my{FINISHED THUS FAR}}
\nc{\FU}{Fr\'echet--Urysohn}
\nc{\gs}{$\gamma$~space}
\nc{\Gab}{\Gamma_{\mathrm{B}}}
\nc{\Omb}{\Omega_{\mathrm{B}}}
\nc{\Ga}{\Gamma}\nc{\Om}{\Omega}
\nc{\smallbinom}[2]{\begin{psmallmatrix} #1\\ #2 \end{psmallmatrix}}
\nc{\bgamma}{\smallbinom{\Om}{\Ga}}
\newcommand{\two}{\{0,1\}}
\nc{\productive}[2]{(#1,\allowbreak #2)^\x}
\nc{\prdct}[1]{#1^\x}
\nc{\Sel}{\mathsf{S}}
\nc{\sset}[2]{\{\,#1 : #2\,\}}
\nc{\smb}[1]{{\!\!\mb{#1}\!\!}}
\nc{\medset}[2]{{\biggl\{\,#1 : #2\,\biggr\}}}
\nc{\smallmedset}[2]{{\bigl\{\,#1 : #2\,\bigr\}}}
\nc{\set}[2]{{\left\{\,#1 : #2\,\right\}}}
\nc{\seq}[2]{{\la\, #1 : #2\,\ra}}
\nc{\eseq}[1]{#1_1, \allowbreak #1_2, \allowbreak\dotsc} 
\nc{\seleseq}[1]{#1_1\in \mathcal{#1}_1, \allowbreak #1_2\in \mathcal{#1}_2, \allowbreak\dotsc}
\nc{\cube}{(\Cantor)^\bbN}
\nc{\Match}{\op{Match}}
\nc{\concat}[1]{\hat{\phantom{a}}\langle #1\rangle}
\nc{\poset}{\mathbb{P}}
\nc{\fn}[1]{{\op{Fn}(#1\times\w,2)}}
\nc{\linadd}{\op{linadd}}
\nc{\nonprod}{\non^\x}
\nc{\alephes}{{\aleph_0}}
\nc{\my}[1]{\marginpar{\textcolor{red}{***}}\textcolor{red}{#1}}
\nc{\later}[1]{{\color{green} #1}}
\nc{\BTs}[1]{{\color{green} #1 (BT)}}
\nc{\Cp}{\op{C}_\mathrm{p}}
\nc{\Bp}{\op{B}_p}
\nc{\Pa}[8]{\bibitem{#1} {#2}, \emph{#3}, {#4} \textbf{#5} ({#6}), {#7}--{#8}.}
\nc{\tPa}[5]{\bibitem{#1} {#2}, \emph{#3}, {#4}, to appear.}
\nc{\sPa}[4]{\bibitem{#1} {#2}, \emph{#3}, {#4}, submitted.}
\nc{\Bc}[9]{\bibitem{#1} {#2}, \emph{#3}, in: \textbf{#4} (#5), #6 #7, #8--#9.}
\nc{\fD}{\mathfrak{D}}
\nc{\fX}{\mathfrak{X}}
\nc{\Onbd}{\Op_{\mathrm{nbd}}} 
\nc{\Omnb}{\Om_{\mathrm{nbd}}} 
\nc{\od}{\mathfrak{od}}
\nc{\Setting}[7]{\xymatrix@R=4pt@C=7pt{#1\ar@{-}[r]&#2\ar@{-}[r]&#3\\&#4\ar@{-}[u]\\
#5\ar@{-}[uu]\ar@{-}[r] & #6\ar@{-}[u]\ar@{-}[r] & #7\ar@{-}[uu]}}
\nc{\mx}[1]{\begin{matrix}#1\end{matrix}}
\nc{\plim}{p\txt{-}\lim}
\nc{\Bgp}{{\Z^\bbN}}
\nc{\Cgp}{{{\Z_2}^\bbN}}
\nc{\Cite}[1]{\textbf{[#1]}}
\nc{\Next}[1]{{#1^+}}
\nc{\cFin}{\mathrm{cF}}
\nc{\scsp}{\text{-scale space}}
\nc{\cfn}{\text{cofinal}\ }
\nc{\Con}{\text{Concentrated}}
\nc{\Lind}{\text{Lindel\"of}\,}
\nc{\con}{\text{-Concentrated}}
\nc{\lind}{\text{-Lindel\"of}\,}
\nc{\ctbl}{\text{countably }\allowbreak}
\nc{\Hur}{\text{Hurewicz}}
\nc{\intvl}[2]{{[#1(#2),\allowbreak #1(#2\!+\!1))}}
\nc{\Bdd}{\mathbf{B}}
\nc{\Dfin}{\mathfrak{D}_\mathrm{fin}}
\nc{\grbl}{{\mbox{\textit{\tiny gp}}}}
\nc{\bbP}{\mathbb{P}}
\nc{\BOfat}{\B_{\Om_{\mathrm{fat}}}}
\nc{\Bgood}{\B_{\mathrm{good}}}
\nc{\compactN}{\cl{\mathbb{N}}}
\nc{\blocks}[2]{\op{cl}_{#2}(#1)}
\nc{\blocksplus}[2]{\op{cl}^+_{#2}(#1)}
\nc{\arx}[1]{\texttt{http://arxiv.org/math/#1}}
\nc{\bq}{\begin{quote}}
\nc{\eq}{\end{quote}}
\nc{\cl}[1]{\overline{#1}}
\nc{\Cl}[2]{\mathrm{cl}_{#1}(#2)}
\nc{\CH}{the Continuum Hypothesis}
\nc{\MA}{Martin's Axiom}
\nc{\Bfat}{\B_\mathrm{fat}}
\nc{\inv}{^{-1}}
\nc{\Cantor}{{\two^\bbN}}
\nc{\bP}{\mathbf{P}}
\nc{\bof}{\op{\fb}}
\nc{\dof}{\op{\fd}}
\nc{\bofF}{\bof(\cF)}
\nc{\sr}[3]{\underset{\mbox{#3}}{\mbox{#1}}}
\nc{\gp}{\binom{\Om}{\Ga}}
\nc{\gpsmall}{\mbox{$\gp$}}
\nc{\gig}{\gimel}
\nc{\gns}{\sone(\Om,\gig)}
\nc{\nsr}[2]{#1}
\nc{\Srg}{{\mathbb{S}}}
\nc{\Srgs}{{\mathbb{S}^*}}
\nc{\NN}{{\bbN^{\bbN}}}
\nc{\ZN}{{\Z^{\bbN}}}
\nc{\NNup}{{\bbN^{\uparrow\bbN}}}
\nc{\NNupb}{{b^{\uparrow\bbN}}}
\nc{\Pof}{\op{P}}
\nc{\PN}{{\Pof(\bbN)}}
\nc{\rothx}[1]{{[#1]^{\mbox{\tiny $\infty$}}}}
\nc{\tx}{{\tilde{x}}}
\nc{\roth}{{[\bbN]^{\mbox{\tiny $\infty$}}}} 
\nc{\roths}{{[b]^{\mbox{\tiny $\infty$}}}} 
\nc{\Fin}{\mathrm{Fin}}
\nc{\ici}{[\bbN]^{ \infty, \infty}}
\nc{\Inc}{{\compactN^{\uparrow\bbN}}}
\nc{\powInc}[1]{{\big(\Inc\big)^{#1}}}
\nc{\powFin}[1]{{\big(\Fin\big)^{#1}}}
\nc{\powPN}[1]{{\big(\PN\big)^{#1}}}
\nc{\NcompactN}{{\compactN^\bbN}}
\nc{\Uarrow}{\smash{\big\uparrow}}
\nc{\LE}{\preccurlyeq}
\nc{\GE}{\succcurlyeq}
\nc{\op}{\operatorname}
\nc{\im}{\op{Im}}
\nc{\Span}{\op{span}}
\nc{\maxfin}{\op{maxfin}}
\nc{\ran}{\op{range}}
\nc{\iso}{\cong}
\nc{\Madd}{{\M}^\star}
\nc{\cI}{\mathcal{I}}
\nc{\cJ}{\mathcal{J}}
\nc{\scrA}{\mathscr{A}}
\nc{\scrB}{\mathscr{B}}
\nc{\scrC}{\mathscr{C}}
\nc{\scrD}{\mathscr{D}}
\nc{\scrF}{\mathscr{F}}
\nc{\scrK}{\mathscr{K}}
\nc{\A}{\D\forall}
\nc{\B}{\mathrm{B}}
\nc{\cB}{\mathcal{B}}
\nc{\cZ}{\mathcal{Z}}
\nc{\bB}{\mathbf{B}}
\nc{\BS}{\mathbf{B}(\mathcal{S})}
\nc{\BF}{\mathbf{B}(\mathcal{F})}
\nc{\BU}{\mathbf{B}(\mathcal{U})}
\nc{\cSp}{\mathcal{S}^+}
\nc{\cFp}{\mathcal{F}^+}
\nc{\cUp}{\mathcal{U}^+}

\nc{\BG}{\B_\Ga}
\nc{\BL}{\B_\Lambda}
\nc{\BT}{\B_\Tau}
\nc{\BTstar}{\B_{\Tau^*}}
\nc{\BO}{\B_\Om}
\nc{\DO}{\cD_\Om}
\nc{\KO}{\cK_\Om}
\nc{\CG}{C_\Ga}
\nc{\CL}{C_\Lambda}
\nc{\CT}{C_\Tau}
\nc{\CTstar}{C_{\Tau^*}}
\nc{\CO}{C_\Om}
\nc{\COgp}{C_{\Om^{\grbl}}}
\nc{\CLgp}{C_{\Lambda^{\grbl}}}
\nc{\BOgp}{\B_{\Om}^{\grbl}}
\nc{\BLgp}{\B_{\Lambda^{\grbl}}}
\nc{\sfC}{\mathsf{C}}
\nc{\sfD}{\mathsf{D}}
\nc{\bD}{\mathbf{D}}
\nc{\Tau}{\mathrm{T}}
\nc{\cA}{\mathcal{A}}
\nc{\cK}{\mathcal{K}}
\nc{\cD}{\mathcal{D}}
\nc{\cF}{\mathcal{F}}
\nc{\cS}{\mathcal{S}}
\nc{\cT}{\mathcal{T}}
\nc{\cG}{\mathcal{G}}
\nc{\cY}{\mathcal{Y}}
\nc{\J}{\mathcal{J}}
\nc{\cL}{\mathcal{L}}
\nc{\cM}{\mathcal{M}}
\nc{\cN}{\mathcal{N}}
\nc{\cH}{\mathcal{H}}
\nc{\cO}{\mathcal{O}}
\nc{\Op}{\mathrm{O}}
\nc{\rmA}{\mathrm{A}}
\nc{\rmF}{\mathrm{F}}
\nc{\rmB}{\mathrm{B}}
\nc{\rmD}{\mathrm{D}}
\nc{\rmP}{\mathrm{P}}
\nc{\cC}{\mathcal{C}}
\nc{\cP}{\mathcal{P}}
\nc{\bbQ}{\mathbb{Q}}
\nc{\bbR}{\mathbb{R}}
\nc{\cU}{\mathcal{U}}
\nc{\cQ}{\mathcal{Q}}
\nc{\Un}{\bigcup}
\nc{\cV}{\mathcal{V}}
\nc{\cW}{\mathcal{W}}
\nc{\Z}{{\mathbb Z}}
\nc{\Impl}{\Rightarrow}
\long\def\forget#1\forgotten{\marginpar{\textcolor{green}{Forgetting...}}}
\nc{\ft}{\mathfrak{t}}
\nc{\fb}{\mathfrak{b}}
\nc{\fc}{\mathfrak{c}}
\nc{\fd}{\mathfrak{d}}
\nc{\fg}{\mathfrak{g}}
\nc{\oo}{\infty}
\nc{\fr}{\mathfrak{r}}
\nc{\fk}{\mathfrak{k}}
\nc{\bidi}{\mathfrak{bidi}}
\nc{\fu}{\mathfrak{u}}
\nc{\fh}{\mathfrak{h}}
\nc{\fp}{\mathfrak{p}}
\nc{\fj}{\mathfrak{j}}
\nc{\fs}{\mathfrak{s}}
\nc{\w}{\omega}
\nc{\x}{\times}
\nc{\Iff}{\Leftrightarrow}

\nc{\nin}{\notin}
\nc{\cat}{\hat{\ }}
\nc{\sub}{\subseteq}
\nc{\spst}{\supseteq}
\nc{\sm}{\setminus}
\nc{\as}{\subseteq^*}
\nc{\les}{\le^*}
\nc{\leinf}{\le^{\infty}}
\nc{\leS}{\le_S}
\nc{\leF}{\le_F}
\nc{\leU}{\le_U}
\nc{\rest}{\restriction}

\nc{\la}{\langle}
\nc{\ra}{\rangle}
\nc{\E}{\exists}
\nc{\dom}{\op{dom}}
\nc{\cov}{\op{cov}}
\nc{\add}{\op{add}}
\nc{\addmen}{\add(\Men{})}
\nc{\cof}{\op{cof}}
\nc{\cf}{\op{cf}}
\nc{\non}{\op{non}}
\nc{\unif}{\op{non}}
\nc{\COV}{\op{COV}}
\nc{\ADD}{\op{ADD}}
\nc{\COF}{\op{COF}}
\nc{\NON}{\op{NON}}
\nc{\impl}{\to}
\nc{\Lp}{\mathcal{L_\p}}
\nc{\Wlog}{without loss of generality}
\newtheorem{thm}{Theorem}[section]
\nc{\bthm}{\begin{thm}} \nc{\ethm}{\end{thm}}
\newtheorem{prop}[thm]{Proposition}
\nc{\bprp}{\begin{prop}} \nc{\eprp}{\end{prop}}
\newtheorem{fact}[thm]{Fact}
\nc{\bfct}{\begin{fact}} \nc{\efct}{\end{fact}}
\newtheorem{prob}[thm]{Problem}
\nc{\bprb}{\begin{prob}} \nc{\eprb}{\end{prob}}
\newtheorem{lem}[thm]{Lemma}
\nc{\blem}{\begin{lem}} \nc{\elem}{\end{lem}}
\newtheorem{app}[thm]{Application}
\nc{\bapp}{\begin{app}} \nc{\eapp}{\end{app}}
\newtheorem{claim}[thm]{Claim}
\nc{\bclm}{\begin{claim}} \nc{\eclm}{\end{claim}}
\newtheorem{cor}[thm]{Corollary}
\nc{\bcor}{\begin{cor}} \nc{\ecor}{\end{cor}}
\newtheorem{conj}[thm]{Conjecture}
\nc{\bcnj}{\begin{conj}} \nc{\ecnj}{\end{conj}}
\theoremstyle{definition}
\newtheorem{defn}[thm]{Definition}
\nc{\bdfn}{\begin{defn}} \nc{\edfn}{\end{defn}}
\newtheorem{obs}[thm]{Observation}
\nc{\bobs}{\begin{obs}} \nc{\eobs}{\end{obs}}
\theoremstyle{remark}
\newtheorem{rem}[thm]{Remark}
\nc{\brem}{\begin{rem}} \nc{\erem}{\end{rem}}
\newtheorem{cnv}[thm]{Convention}
\nc{\bcnv}{\begin{cnv}} \nc{\ecnv}{\end{cnv}}
\newtheorem{exam}[thm]{Example}
\nc{\bexm}{\begin{exam}} \nc{\eexm}{\end{exam}}
\nc{\bpf}{\begin{proof}} \nc{\epf}{\end{proof}
}
\nc{\be}{\begin{enumerate}}
\nc{\ee}{\end{enumerate}}
\nc{\bi}{\begin{itemize}}
\nc{\bimy}{\my{\begin{itemize}}
\nc{\eimy}{\end{itemize}}}
\nc{\itm}{\item}
\nc{\ei}{\end{itemize}}
\nc{\Subsection}[1]{\goodbreak\subsection*{#1}}

\nc{\sone}{\mathsf{S}_1}
\nc{\sfin}{\mathsf{S}_\mathrm{fin}}
\nc{\ufin}{\mathsf{U}_\mathrm{fin}}
\nc{\Split}{\mathsf{Split}}

\nc{\gone}{\mathsf{G}_1}    \nc{\gfin}{\mathsf{G}_\mathrm{fin}}

\nc{\men}{\sfin(\Op,\Op)}
\nc{\sch}{\ufin(\cO,\Omega)}
\nc{\rothb}{\text{Rothberger}}
\nc{\pmen}{\sfin(\Omega,\Omega)}
\nc{\Rothb}{\sone(\Op,\Op)}

\nc{\prothb}{\sone(\Omega,\Omega)}
\nc{\tU}{{\tilde{U}}}
\nc{\tF}{{\tilde{F}}}
\nc{\tY}{{\tilde{Y}}}
\nc{\tX}{{\tsup[1]{X}}}
\nc{\dtX}{{\tsup[2]{X}}}
\nc{\dt}[1]{{\tsup[2]{#1}}}
\nc{\td}{{\tilde{d}}}
\nc{\tz}{{\tilde{z}}}
\nc{\cfd}{\cf(\fd)}

\nc{\msep}{\sfin(\cD,\cD)}
\nc{\rsep}{\sone(\cD,\cD)}
\nc{\cft}{\sfin(\Omega_{\mathbf{0}},\Omega_{\mathbf{0}})}
\nc{\scft}{\sone(\Omega_{\mathbf{0}},\Omega_{\mathbf{0}})}

\nc{\Umen}{U\text{-Menger}}
\nc{\hur}{\ufin(\Op,\Gamma)}
\nc{\tUmen}{\tU\text{-Menger}}

\nc{\Men}{\text{Menger}}
\nc{\Sch}{\text{Scheepers}}

\nc{\aspst}{\prescript{*}{}{\spst}\ }
\nc{\eqs}{=^*}
\nc{\ctblOm}{\Omega_{\mathrm{ctbl}}}
\nc{\GNga}{{\smallbinom{\Om}{\Ga}}}
\nc{\ctblga}{\smallbinom{\ctblOm}{\Ga}}

\nc{\nadd}{\cN_{\mathrm{add}}}
\nc{\ball}{\mathrm{B}}
\nc{\cOunif}{\cO^{\textrm{unif}}}

\nc{\sep}{
\vspace{2cm}
\noindent
\begin{minipage}{\textwidth}
	\textcolor{red}{\rule{\textwidth}{1pt}}
\end{minipage}
}

\nc{\nonNadd}{\non(\nadd)}
\nc{\borga}{\Ga_\mathrm{Bor}}
\nc{\pick}{x}
\nc{\gen}{y}
\nc{\nullzind}{\sone(\{\Op_n^{\mathsf{unif}}\}_{n\in\bbN},\Ga)}
\nc{\nullzindf}[1]{\sone(\{\Op_{#1}^{\mathsf{unif}}\}_{n\in\bbN},\Ga)}

\nc{\myb}[1]{\marginpar{\textcolor{blue}{***}}\textcolor{blue}{#1}}